\documentclass[10pt,reqno]{amsart}
\usepackage{amsmath,multicol}
\usepackage{amsthm}
\usepackage{amssymb}
\usepackage{amsfonts}
\usepackage{subcaption}
\usepackage{graphicx}
\usepackage{latexsym}
\usepackage{cite, url, xcolor}
\usepackage{stmaryrd}
\usepackage{hyperref}
\usepackage{enumitem}
    \setenumerate{itemsep=5pt} 
    \setitemize{itemsep=5pt} 

\usepackage[font=Large]{caption}

\usepackage[font=footnotesize,labelfont=bf]{caption}

\usepackage[sc]{mathpazo}
\usepackage[T1]{fontenc}

\usepackage{tikz,amsmath,amsfonts}
\usetikzlibrary{positioning}
\usetikzlibrary{decorations.pathreplacing}

\newcommand{\la}{\lambda}

\newcommand{\RR}{\mathbb{R}}
\newcommand{\CC}{\mathbb{C}}

\newcommand{\cP}{\mathcal{P}}

\newtheorem{Theorem}{Theorem}

 \newtheorem{Proposition}[Theorem]{Proposition}
 \newtheorem{Example}[Theorem]{Example}

\allowdisplaybreaks
\begin{document}
\title[B-Splines, Positivity, and Complete Homogeneous Symmetric Polynomials]
{Weighted Means of B-Splines, Positivity of Divided Differences, and Complete Homogeneous Symmetric Polynomials}

\author[A.~B\"{o}ttcher]{Albrecht B\"{o}ttcher}
\address{Fakult\"{a}t f\"{u}r Mathematik, Technische Universit\"at Chemnitz, 09107 Chemnitz, Germany}
\email{aboettch@mathematik.tu-chemnitz.de}
\urladdr{\url{https://www-user.tu-chemnitz.de/~aboettch/}}

\author[S.R.~Garcia]{Stephan Ramon Garcia}
\address{Department of Mathematics, Pomona College, 610 N. College Ave., Claremont, CA 91711}
\email{stephan.garcia@pomona.edu}
\urladdr{\url{http://pages.pomona.edu/~sg064747}}

\author[M.~Omar]{Mohamed Omar}
\address{Department of Mathematics, Harvey Mudd College, 301 Platt Blvd., Claremont, CA 91711}
\email{omar@g.hmc.edu}
\urladdr{\url{www.math.hmc.edu/~omar}}

\author[C.~O'Neill]{Christopher O'Neill}
\address{Mathematics Department, San Diego State University, San Diego, CA 92182}
\email{cdoneill@sdsu.edu}
\urladdr{\url{https://cdoneill.sdsu.edu/}}

\thanks{SRG partially supported by NSF grant DMS-1800123.}

\begin{abstract}
We employ the fact certain divided differences can be written as weighted means of B-splines
and hence are positive. These divided differences include the complete homogeneous symmetric polynomials
of even degree $2p$, the positivity of which is a classical result by D. B. Hunter. We extend Hunter's result
to  complete homogeneous symmetric polynomials of fractional degree, which are defined via Jacobi's
bialternant formula. We show in particular that these polynomials have positive real part
for real degrees $\mu$ with $|\mu-2p|< 1/2$. We also prove a positivity criterion for linear
combinations of the classical complete homogeneous symmetric polynomials and a sufficient criterion
for the positivity of linear combinations of products of such polynomials.
\end{abstract}

\maketitle

\section{Introduction}
\label{absec1}

\noindent
A formula by Peano states that if $f: [a_1,a_n] \to \RR$ is an $n-1$ times continuously differentiable function, then
\begin{equation}
f[a_1, \ldots,a_n]=\frac{1}{(n-1)!}\int_{a_1}^{a_n} f^{(n-1)}(x) F(x;a_1, \ldots,a_n)\,dx. \label{bo1}
\end{equation}
See, for instance, \cite[Chap. III, Sec. 3.7]{Davis} or~\cite{CS}. Here $a_1 < a_2 < \cdots < a_n$ are real numbers, $f[a_1, \ldots,a_n]$
denotes the $n$th divided difference of $f$, which may be written as
\[f[a_1, \ldots,a_n]=\sum_{j=1}^n \frac{f(a_j)}{\prod_{k\neq j}(a_j-a_k)},\]
and $F(x;a_1,\ldots,a_n)$ is the Curry--Schoenberg B-spline introduced in~\cite{CS}, one representation of which is
\begin{equation}
F(x;a_1,a_2,\ldots,a_n) = \frac{n-1}{2} \sum_{j=1}^n \frac{|a_j-x|(a_j-x)^{n-3}}{\prod_{k \neq j}(a_j-a_k)}. \label{bo2}
\end{equation}
The function $F(x;a_1,\ldots,a_n)$ is positive on $(a_1,a_n)$, and hence if $f^{(n-1)}$ is nonnegative and not identically zero
on $\RR$, then~(\ref{bo1}) implies that $f[a_1,\ldots,a_n] >0$ whenever $a_1 < a_2 < \cdots < a_n$. Taking $f(x)=x^{p+n-1}$ with
a nonnegative integer $p$, we obtain from~(\ref{bo1}) that
\begin{equation}
f[a_1, \ldots,a_n]=\binom{p+n-1}{n-1}\int_\RR x^{p} F(x;a_1, \ldots,a_n)\,dx, \label{bo3}
\end{equation}
and it is well known
that in this case $f[a_1, \ldots,a_n]$ is just the complete homogeneous symmetric polynomial
\begin{equation*}
h_p(a_1,a_2,\ldots,a_n) = \sum_{1 \leq j_1 \leq j_2 \leq \cdots \leq j_p \leq n} a_{j_1} a_{j_2} \cdots a_{j_p},
\end{equation*}
with the convention $h_0(a_1, a_2,\ldots, a_n)=1$; see, e.g.,~\cite[Theorem 1.2.1]{Phil}. Thus, if $p$ is an even positive integer then $h_p(a_1,a_2,\ldots,a_n) >0$
for arbitrary pairwise different real numbers $a_1, \ldots,a_n$. This is a classical result by Hunter~\cite{Hunter}. He proved it in a completely different way.
The easy argument employed above is essentially from~\cite{GOO}.

The purpose of this paper is to extend Hunter's result to more general functions, in particular to complete homogeneous symmetric
polynomials of fractional degree and to sums of products of the classical homogeneous symmetric polynomials. Our approach is
based on the preceding argument and a determinantal representation of B-splines.

\section{Main Results}
\label{bosec1}

\noindent
Jacobi's bialternant formula says that for each positive integer $z$ we have
\begin{equation} \label{abJac}
h_z(a_1,a_2,\ldots,a_n)  V(a_1,a_2,\ldots,a_n) = \det
\begin{bmatrix}
1 & a_1 & a_1^2 & \cdots & a_1^{n-2} & a_1^{z+n-1} \\
1 & a_2 & a_2^2 & \cdots & a_2^{n-2} & a_2^{z+n-1} \\
\vdots & \vdots & \vdots & \ddots & \vdots & \vdots \\
1 & a_n & a_n^2 & \cdots & a_n^{n-2} & a_n^{z+n-1} \\
\end{bmatrix},
\end{equation}
where $h_z(a_1,a_2,\ldots,a_n)$ is the complete homogeneous symmetric polynomial introduced above
and
\[
V(a_1,a_2,\ldots,a_n)=
\det \begin{bmatrix}
1 & a_1 & a_1^2 & \cdots & a_1^{n-1} \\
1 & a_2 & a_2^2 & \cdots & a_2^{n-1} \\
\vdots & \vdots & \vdots & \ddots & \vdots \\
1 & a_n & a_n^2 & \cdots & a_n^{n-1} \\
\end{bmatrix}=\prod_{1 \leq i < j \leq n} (a_j-a_i)
\]
is the $n \times n$ Vandermonde determinant; see, e.g., \cite{ec2}.

For the moment we assume that $a_1, a_2, \ldots,a_n$
are pairwise distinct real numbers and that none of them is zero.
Putting $a_j^z=e^{z \log a_j}$ with a choice of the branch of the logarithm that is defined on $\RR\setminus\{0\}$, the
right-hand side makes sense for every $z \in \CC$.
Thus,
it is natural to take~(\ref{abJac}) as the definition of $h_z(a_1, a_2,\ldots, a_n)$ for $z \in \CC$,
that is, of {\em fractional degree complete homogeneous symmetric polynomials}.
For example, in the case of three variables we obtain
\begin{equation}\label{abh3}
h_z(a,b,c)
=\frac{a^{z+2}(b-c)   +b^{z+2}(c-a) + c^{z+2}(a-b) }{ (a-b) (a-c) (b-c)}.
\end{equation}
As said, for positive integers $z$,  these are the usual complete homogeneous symmetric polynomials.
For other choices of $z$ we get, for instance,
$h_{-1}(a,b,c) = h_{-2}(a,b,c) = 0$, and
\begin{align*}
h_{\frac{1}{2}}(a,b,c)
&= \frac{ a^{\frac{5}{2}}(b-c)   +b^{\frac{5}{2}}(c-a) + c^{\frac{5}{2}}(a-b) }{(a-b) (a-c) (b-c)} ,\\
h_{-\frac{1}{2}}(a,b,c)
&= \frac{\sqrt{a} \sqrt{b}+\sqrt{a} \sqrt{c}+\sqrt{b} \sqrt{c}}{(\sqrt{a}+\sqrt{b}) (\sqrt{a}+\sqrt{c}) (\sqrt{b}+\sqrt{c})} ,\\
h_{-\frac{3}{2}}(a,b,c)
&=-\frac{1}{(\sqrt{a}+\sqrt{b}) (\sqrt{a}+\sqrt{c}) (\sqrt{b}+\sqrt{c})} ,\\
h_{-\frac{5}{2}}(a,b,c)
&= \frac{\frac{a-b}{\sqrt{c}}+\frac{b-c}{\sqrt{a}}+\frac{c-a}{\sqrt{b}}}{(a-b) (a-c) (b-c)} ,\\
h_{-3}(a,b,c)
&=\frac{1}{abc} ,\\
h_{-4}(a,b,c)
&= \frac{a b+a c+b c}{a^2 b^2 c^2},\\
h_i(a,b,c)& = \frac{a^2 e^{i\log a}(b-c)   +b^2 e^{i \log b}(c-a) + c^2 e^{i \log c} (a-b) }{ (a-b) (a-c) (b-c)},\quad i =\sqrt{-1},
\end{align*}
each of which is a symmetric function of $a,b,c$. These examples are already in~\cite{GOO}, where an unexpected connection
between complete homogeneous symmetric polynomials and the statistical properties of factorization lengths in numerical semigroups
was investigated.

If ${\rm Re}\,(z+n-1) >0$, we put $0^{z+n-1}:=0$. Thus, under the assumption that ${\rm Re}\,z > -1$, we may use (\ref{abJac}) to define $h_z(a_1,a_2, \ldots,a_n)$
for $n \ge 2$ also in the situation where (exactly) one of the numbers $a_1, a_2,\ldots, a_n$ is zero. The following theorem provides us with
an alternative representation of complete homogeneous symmetric polynomials.

\begin{Theorem}\label{abthm1}
Let $n \ge 2$, let $a_1,a_2,\ldots,a_n \in \mathbb{R}$ with $a_1 < a_2 < \cdots < a_n$, and let $F(x;a_1, \ldots,a_n)$
be the function~\eqref{bo2}.
If $z \in \CC$ and ${\rm Re}\,z >-1$, then
$x^z  F(x;a_1,a_2,\ldots,a_n)$ is absolutely integrable and
\begin{equation}\label{abhz}
 h_z (a_1,a_2,\ldots,a_n)=\binom{z+n-1}{n-1} \int_{\RR} x^z  F(x;a_1,a_2,\ldots,a_n)\,dx.
\end{equation}
\end{Theorem}

The use
of such a theorem in connection with fractional degree complete homogeneous symmetric polynomials
was first indicated in~\cite{GOO}. There it was shown that $F(x; a_1, \ldots, a_n)$ is a probability density supported
in $[a_1,a_n]$ which is piecewise polynomial of degree $n-2$ and which is $n-3$ times continuously differentiable.
In~\cite{GOO}, the function arose in the formula
\[
\lim_{m \to \infty}\frac{|\{\ell \in {\sf L}[\![m]\!]: \ell \in [\alpha m, \beta m]\}|}{|{\sf L}[\![m]\!]|}=
\int_\alpha^\beta F(x;1/m_n, \ldots,1/m_1)\, dx,\]
where ${\sf L}[\![m]\!]$ is the multiset of lengths $\ell = x_1+\cdots+x_n$ of possible decompositions $m=x_1m_1+\cdots+x_nm_n$
with nonnegative integers $x_j$ for given positive integers $m_1<\cdots< m_n$ satisfying ${\rm gcd}(m_1, \ldots, m_n)=1$. (This is related
to the coin problem of Frobenius.)

Once Grigori Olshanski saw a preliminary version of this paper, he kindly informed us that $F(x;a_1, \ldots,a_n)$
is nothing but the B-spline introduced by Curry and Schoenberg in~\cite{CS}. Thanks to this observation we were released
from our effort devoted to proving positivity, support in $[a_1,a_n]$, and unimodality of the function $F(x;a_1, \ldots,a_n)$
since these turned out to be well-known properties of B-splines. As Olshanski pointed out, with $x_+:=\max(x,0)$ we
have $|x|=2x_+-x$ and $x_+x^{n-3}=x_+^{n-2}$,
hence the sum
in~(\ref{bo2}) equals
\[\frac{n-1}{2}\sum_{j=1}^n \frac{2(a_j-x)_+^{n-2}}{\prod_{k\neq j}(a_j-a_k)}-\frac{n-1}{2}\sum_{j=1}^n \frac{(a_j-x)^{n-2}}{\prod_{k\neq j}(a_j-a_k)},\]
and since the second sum is just (\ref{bo2}) for $x < a_1$, which is known to be zero, we get
\[F(x;a_1, \ldots,a_n)=(n-1)\sum_{j=1}^n \frac{(a_j-x)_+^{n-2}}{\prod_{k\neq j}(a_j-a_k)},\]
which is exactly the formula given in~\cite{CS} and in~\cite{Ol}.

Curry and Schoenberg proved in particular that $F(x;a_1, \ldots,a_n)$
is a probability density supported in $[a_1,a_n]$ and that the $k$th
derivative ($k=0,1,\ldots, n-3$) of the function has exactly $k$ simple zeros in $(a_1,a_n)$.
They also proved the remarkable geometric interpretation of $F(x;a_1, \ldots,a_n)$ as the linear density function obtained
by projecting orthogonally onto the $x$-axis
the volume of an $(n-1)$-dimensional simplex of volume $1$, so located
that its $n$ vertices project orthogonally into the points $a_1, \ldots, a_n$ of
the $x$-axis. That an interpretation of this type might be true was independently communicated
to us by Terence Tao.

For more on splines we refer
to the monographs~\cite{Boor} and~\cite{Phil}. We here only note that there are different normalizations of B-splines: those
of Curry and Schoenberg are normalized so that their integral is $1$ whereas the B-Splines of de Boor are normalized so that
a certain collection of them (over shifted intervals) sums to $1$.

So far we assumed that $a_1 < a_2 < \cdots < a_n$. By appropriate limit passages or by constructing B-splines via
recursion formulas and making thorough use of the convention $0/0:=0$ (called ``the useful maxim'' on page~117 of~\cite{Boor}),
one may extend the definition of $F(x; a_1, \ldots,a_n)$
to arbitrary $a_1 \le a_2 \le \cdots \le a_n$ under the mere assumption that $a_1 < a_n$. The resulting
functions are still positive piecewise-polynomial probability densities supported in $[a_1,a_n]$, and
only the smoothness is lowered to some $n-r < n-3$.

For a positive integer $p$ the polynomial $h_p(a_1, \ldots, a_n)$
is well-defined without the assumption that the $a_j$ be pairwise distinct. Again by appropriate limit passages in Vandermonde-like
determinants (leading to so-called confluent Vander\-monde-like determinants), one may also
define $h_z(a_1, \ldots,a_n)$ for ${\rm Re}\, z > -1$ under the sole requirement that among $a_1, \ldots, a_n$ there are at least
two different numbers. Finally, for $a \neq 0$, the natural definition of $h_z(a, \ldots,a)$ respecting continuity is
\begin{equation}
h_z(a, \ldots, a)=\binom{z+n-1}{n-1}a^z=\frac{(z+n-1)\cdots(z+1)}{(n-1)!}a^z.\label{aaaa}
\end{equation}

These limit passages give Theorem \ref{abthm1} under the only assumption that
\begin{equation} \label{aaaaa}
a_1 \le a_2 \le \ldots \le a_n \;\:\mbox{and}\;\: a_1 < a_n.
\end{equation}

The classical result by D. B. Hunter \cite{Hunter} mentioned states that if $p$ is a nonnegative integer, then $h_{2p}(a_1, \ldots, a_n) >0$
for all $(a_1, \ldots,a_n) \in \RR^n\setminus\{(0,\ldots,0)\}$. See~\cite{Tao2} for more results on this topic.
We will prove the following generalization of Hunter's result.

\begin{Theorem}\label{abthm2}
Choose the branch of the complex logarithm that is analytic on $\CC$ cut along the negative imaginary axis and takes the value $0$ at $1$.
Let $\mu > -1$ be a real number  and suppose
$(a_1, \ldots,a_n) \in \mathbb{R}^n \setminus \{(0,\ldots,0)\}$.

\medskip
{\rm (a)} If  $|\mu-2p| < \tfrac{1}{2}$ for some nonnegative integer $p$, then
${\rm Re}\,h_{\mu}(a_1, \ldots, a_n) >0$.

\medskip
{\rm (b)} If $|\mu-(2p-1)|<\tfrac{1}{2}$ for some integer $p \ge 0$, then ${\rm Re}\,h_\mu(a_1, \ldots,a_n)>0$ for
$(a_1, \ldots,a_n) \in [0,\infty)^n$ and ${\rm Re}\,h_\mu(a_1, \ldots,a_n)<0$ for
$(a_1, \ldots,a_n) \in (-\infty,0]^n$.

\medskip
{\rm (c)} If $|\mu-p|=\tfrac{1}{2}$ for some nonnegative integer $p$, then ${\rm Re}\,h_\mu(a_1, \ldots,a_n)\ge 0$, and we have
${\rm Re}\,h_\mu(a_1, \ldots,a_n)=0$ for
$(a_1, \ldots,a_n) \in (-\infty,0]^n$.
\end{Theorem}

Note that that the cases $|\mu-2p|<\tfrac{1}{2}$, $|\mu-(2p+1)|<\tfrac{1}{2}$, and $|\mu-p|=\tfrac{1}{2}$ are equivalent to the cases
$\cos(\mu\pi)>0$, $\cos(\mu\pi)<0$, and $\cos(\mu\pi)=0$, respectively. Section~\ref{absec3} contains some more results related to
Theorem~\ref{abthm2}.
Theorems~\ref{abthm1} and~\ref{abthm2} complement recent work of Terence~Tao~\cite{Tao2} concerning different ways of proving the positivity
of even degree complete homogeneous symmetric polynomials.  We emphasize that the polynomials considered here are polynomials of fractional degree
and that they should be distinguished from the symmetric functions in a fractional number of variables introduced implicitly in~\cite{Tao3}.

We now turn to combinations of the classical complete homogeneous symmetric polynomials. The following result is about linear combinations.

\begin{Theorem} \label{abthm3}
Let $H(a_1, a_2, \ldots, a_n)=\sum_{j=0}^m c_j h_j(a_1, a_2, \ldots, a_n)$ with real coefficients $c_j$ and let $-\infty \le r < s \le \infty$. Then
\[H(a_1, \ldots, a_n) >0\;\:\mbox{for all}\;\: (a_1,\ldots,a_n) \in (r,s)^n\setminus\{(0,\ldots,0)\}\]
if and only if
$H(a,a, \ldots,a) >0$ for all  $a \in (r,s)\setminus\{0\}$.
\end{Theorem}

Here is a sufficient condition for the positivity of combinations involving products. We confine ourselves to the case
of at most two factors. The extension to more than two factors is obvious.

\begin{Theorem} \label{abthm4}
Let
\[H(a_1, \ldots, a_n)=\sum_{j,k=1}^m c_{jk}h_j(a_1, \ldots,a_n)h_k(a_1, \ldots, a_n)\]
with real coefficients $c_{jk}$ and let $-\infty \le r < s \le \infty$. Put
\[\cP(x,y)=\sum_{j,k=1}^m c_{jk} h_j(x,x, \ldots,x) h_k(y,y, \ldots, y).\]
If $\cP(x,y) \ge 0$ for $(x,y) \in (r,s)^2 \setminus\{(0,0)\}$ and $H(a,a, \ldots, a)>0$
for $a \in (r,s) \setminus\{0\}$, then $H(a_1, \ldots, a_n) >0$ for $(a_1, \ldots, a_n)$ in
$(r,s)^n \setminus\{(0,\ldots,0)\}$.
\end{Theorem}

We remark that Theorem \ref{abthm4} is subtler than it might appear at the first glance. Consider, for example,
\[H(a_1,\ldots,a_n)=2\alpha h_2(a_1, \ldots,a_n) h_4(a_1, \ldots,a_n)-3\beta h_2^2(a_1, \ldots,a_n)+2,\]
and let us omit the arguments, that is, let us simply write
\[H=2\alpha h_2h_4-3\beta h_2^2+2.\]
Recall that $h_0(a_1, \ldots, a_n)=1$, so that $2$ may be interpreted as $2h_0^2$.
By (\ref{aaaa}), the polynomial $\cP(x,y)$ is
\[2\alpha \binom{2+n-1}{n-1}\binom{4+n-1}{n-1} x^2 y^4-3 \beta \binom{2+n-1}{n-1}^2 x^2y^2+2,\]
and let us choose $\alpha$ and $\beta$ so that this becomes
\[\cP(x,y)=2 x^2 y^4 - 3x^2 y^2 +2.\]
Since $\cP(x,1)=-x^2+2 < 0$ for $x >\sqrt{2}$, Theorem \ref{abthm4} does not give anything for $(r,s)=(-\infty,\infty)$.
However, we may write
\[H=\alpha h_4 h_2 + \alpha h_2 h_4 -3 \beta h_2^2+2,\]
and now, with the same choice of $\alpha$ and $\beta$ as above, we obtain
\[\cP(x,y)=x^4y^2+x^2y^4-3x^2y^2+2.\]
This is $1$ plus the famous Motzkin polynomial. (The Motzkin polynomial introduced in~\cite{Motz}
was the first explicit example of a nonnegative polynomial that is not a sum of squares of polynomials.
See~\cite{Ben} for a recent survey. Note that nonnegativity is simple: we have
\[x^2y^2 =\sqrt[3]{x^4y^2 \cdot x^2y^4\cdot 1} \le \frac{1}{3}(x^4y^2+x^2y^4+1)\]
by the arithmetic-geometric mean inequality.)
Hence $\cP(x,y) \ge 1$ on all of $\RR^2$.
As also $H(a,a, \ldots,a)=\cP(a,a)\ge 1$ for all $a$, we can now invoke Theorem~\ref{abthm4}
to conclude that $H(a_1, \ldots,a_n) >0$ for all $(a_1, \ldots,a_n) \in \RR^n$. One is tempted to draw
this conclusion from inserting $u=h_2$ and $v=h_4$ in the inequality
\[g(u,v)=2\alpha u v -3 \beta u^2+2 >0 \;\:\mbox{for}\;\: (u,v)\in [0,\infty)^2,\]
but this inequality is not true because $g(u,1) \to -\infty$ as $u \to \infty$.

Theorems \ref{abthm1} to \ref{abthm4} will be proved in Sections~\ref{absec2} and~\ref{bosec5}. In Section~\ref{absec4}
we establish expressions for $h_z(a_1,\ldots,a_n)$ in terms of Schur polynomials in the cases where $z$ is a negative integer or a positive
rational number.

\section{Proof of Theorem \ref{abthm1}}
\label{absec2}

\noindent
We first rewrite $F(x;a_1,a_2, \ldots,a_n)$ in terms of determinants.
Let $a_1,a_2,\ldots,a_n \in \RR$ with $a_1<a_2<\cdots<a_n$ and let $F(x;a_1, a_2, \ldots,a_n)$ be defined by \eqref{bo2}.
In what follows, $V(a_1,\ldots, \widehat{a}_j, \ldots, a_n)$ denotes the $(n-1) \times (n-1)$
Vandermonde determinant obtained from $V(a_1, a_2, \ldots,a_n)$ by removing $a_j$.
Then
\begin{eqnarray*}
& & F(x;a_1,a_2,\ldots,a_n)\\
& & \quad =\frac{n-1}{2} \sum_{j=1}^n \frac{|a_j-x|(a_j-x)^{n-3}}{\prod_{k \neq j}(a_j-a_k)} \\
& & \quad=\frac{n-1}{2} \sum_{j=1}^n \frac{|a_j-x|(a_j-x)^{n-3}}{\prod_{1 \leq k < j}(a_j-a_k)\prod_{j < k \leq n}(a_j-a_k)} \\
& & \quad=\frac{(-1)^{n-j}(n-1)}{2} \sum_{j=1}^n \frac{|a_j-x|(a_j-x)^{n-3}}{\prod_{1 \leq k < j}(a_j-a_k)\prod_{j < k \leq n}(a_k-a_j)} \\
& & \quad= \frac{n-1}{2} \sum_{j=1}^n \frac{ V(a_1,\ldots,\widehat{a_j},\ldots,a_n)}{ V(a_1,a_2,\ldots,a_n)}
(-1)^{n-j} |a_j-x|(a_j-x)^{n-3} \\
& &\quad = (n-1) \frac{\sum_{j=1}^n (-1)^{n+j} V(a_1,\ldots,\widehat{a_j},\ldots,a_n) \cdot
 |a_j-x|(a_j-x)^{n-3}}{2 V(a_1,a_2,\ldots,a_n)} \\
\end{eqnarray*}
and hence
\begin{eqnarray}
& & \!\!\!\!\!F(x;a_1,a_2,\ldots,a_n) \nonumber\\
& & \!\!\!\!\!=  \frac{n-1}{2  V(a_1,a_2,\ldots,a_n)} \det
 \begin{bmatrix}
1 & {a_1} & a_1^2 & \cdots &  a_1^{n-2}  & |a_1-x|(a_1-x)^{n-3} \\
1 & {a_2} & a_2^2 & \cdots & a_2^{n-2} & |a_2-x|(a_2-x)^{n-3} \\
\vdots & \vdots & \vdots & \ddots & \vdots  & \vdots \\
1 & {a_n} & a_n^2 & \cdots &  a_n^{n-2} &  |a_n-x|(a_n-x)^{n-3} \\
\end{bmatrix}. \label{abeq:det}
\end{eqnarray}

We now prove (\ref{abhz}), that is, the equality
\[g_z(a_1, \ldots, a_n)=\binom{z+n-1}{n-1}\int_\RR x^z f(a_1, \ldots,a_n)\,dx\]
with
\begin{equation} \label{aabg}
g_z(a_1, \ldots, a_n)= \det \begin{bmatrix}
1 & a_1 & a_1^2 & \cdots & a_1^{n-2} & a_1^{z+n-1} \\
1 & a_2 & a_2^2 & \cdots & a_2^{n-2} & a_2^{z+n-1} \\
\vdots & \vdots & \vdots & \ddots & \vdots & \vdots \\
1 & a_n & a_n^2 & \cdots & a_n^{n-2} & a_k^{z+n-1} \\
\end{bmatrix}
\end{equation}
and
\begin{equation}
f(a_1, \ldots,a_n)= \frac{n-1}{2} \det
 \begin{bmatrix}
1 & {a_1} & a_1^2 & \cdots &  a_1^{n-2}  & |a_1-x|(a_1-x)^{n-3} \\
1 & {a_2} & a_2^2 & \cdots & a_2^{n-2} & |a_2-x|(a_2-x)^{n-3} \\
\vdots & \vdots & \vdots & \ddots & \vdots  & \vdots \\
1 & {a_n} & a_n^2 & \cdots &  a_n^{n-2} &  |a_n-x|(a_n-x)^{n-3} \\
\end{bmatrix}.
\label{aabf}
\end{equation}
We may assume that $a_j \neq 0$ for all $j$ because
both~(\ref{aabg}) and~(\ref{aabf}) depend continuously on $a_1, \ldots, a_j$. Multiplying~(\ref{aabf})
by $x^z$ and integrating the result amounts to replacing the $j$th entry of the last column by
\begin{eqnarray*}
& & \int_{a_1}^{a_n} x^z|a_j-x|(a_j-x)^{n-3}\,dx\\
& & \quad =\int_{a_1}^{a_j}x^z(a_j-x)^{n-2}\,dx - \int_{a_j}^{a_n}x^z(a_j-x)^{n-2}\,dx\\
& & \quad = \left(\int_0^{a_j}-\int_0^{a_1}-\int_0^{a_n}+\int_0^{a_j}\right)x^z (a_j-x)^{n-2}\,dx\\
& & \quad = 2\int_0^{a_j}x^z(a_j-x)^{n-2}\,dx -\int_0^{a_1}x^z(a_j-x)^{n-2}\,dx
-\int_0^{a_n}x^z(a_j-x)^{n-2}\,dx\\
& & \quad =: 2I_1-I_2^j-I_3^j.
\end{eqnarray*}
We have
\[
I_2^j  =  \sum_{k=0}^{n-2}\int_0^{a_1} x^z \binom{n-2}{k}a_j^k (-1)^{n-2-k}x^{n-2-k}\,dx
 =  \sum_{k=0}^{n-2} c_k(z)a_j^k
\]
and, analogously, $I_3^j=\sum_{k=0}^{n-2}d_k(z)a_j^k.$ It follows that the columns ${\rm col}(I_2^j)_{j=1}^n$
and ${\rm col}(I_3^j)_{j=1}^n$ are linear combinations of the first $n-1$ columns of the determinant~(\ref{aabf}).
Consequently, the $j$th entry of the multiplied and integrated determinant may simply replaced by $2I_1$.
We finally have
\begin{eqnarray*}
2I_1 & = & 2\int_0^{a_j}x^z(a_j-x)^{n-2}\,dx = 2 a_j^{z+n-1}\int_0^1 t^z (1-t)^{n-2}\,dt\\
& = & 2 a_j^{z+n-1}\frac{\Gamma(z+1)\Gamma(n-1)}{\Gamma(z+n)} =
2 a_j^{z+n-1} \frac{\Gamma(z+1)(n-2)!}{(z+n-1) \cdots (z+1)\Gamma(z+1)}\\
& = & 2 a_j^{z+n-1} \binom{z+n-1}{n-1}^{-1}\frac{1}{n-1},
\end{eqnarray*}
which is the asserted equality. $\;\:\square$

\section{Proof of and More Results Around Theorem \ref{abthm2}}
\label{absec3}

\noindent
{\em Proof of Theorem~\ref{abthm2}.}
Since $h_\mu(a)=h_\mu(0,0,a)$ and $h_\mu(a,b)=h_\mu(0,a,b)$, we may restrict ourselves to $n \ge 3$.
Suppose first that all $a_j$ are equal to $a \neq 0$. Then $a^\mu >0$ for $a>0$,
and for $a <0$ we have
\[a^\mu=e^{\mu \log a}=e^{\mu (\log|a|+i\arg a)} =e^{\mu(\log |a|+i \pi)}=|a|^\mu \cos(\mu \pi)+i|a|^\mu \sin (\mu \pi).\]
Consequently, (\ref{aaaa}) implies all assertions of the theorem. As~(\ref{abhz}) remains true if~(\ref{aaaaa}) holds, we obtain that
\[h_\mu(a_1, \ldots,a_n)=\frac{(\mu+n-1)\cdots (\mu+1)}{(n-1)!}\int_\RR x^\mu F(x;a_1, \ldots,a_n)\,dx.\]
With $F(x;a_1, \ldots,a_n)$ abbreviated to $F(x)$, it follows that ${\rm Re}\,h_\mu(a_1,\ldots,a_n)$ is a positive
constant times
\begin{eqnarray}
& &  {\rm Re}\,\left(\int_{-\infty}^0 e^{i\mu\pi}|x|^\mu F(x)\,dx +\int_0^\infty |x|^\mu F(x)\,dx\right)\nonumber\\
& & \quad = \cos(\mu\pi)\int_{-\infty}^0 |x|^\mu F(x)\,dx+\int_0^\infty |x|^\mu F(x)\,dx. \label{abfin}
\end{eqnarray}
If $\cos(\mu\pi) >0$, then (\ref{abfin}) is greater than or equal to $\cos(\mu\pi) \int_\RR|x|^\mu F(x)\,dx$,
and this is strictly greater than zero because $F(x) > 0$ on some open interval. Let $\cos(\mu\pi) <0$. If $a_1 \ge 0$, then ~(\ref{abfin}) equals
$\int_\RR |x|^\mu F(x)\,dx$, which is strictly positive because $F(x)$ is strictly positive on some open interval, and if $a_n \le 0$,
then (\ref{abfin}) is $\cos(\mu\pi)\int_\RR |x|^\mu F(x)\,dx$, which now is strictly negative. Finally, if $\cos(\pi\mu)=0$,
then~(\ref{abfin}) equals $\int_0^\infty |x|^\mu F(x)\,dx$. This is always nonnegative and this vanishes if $a_n \le 0$.
$\;\:\square$

\medskip
Hunter~\cite{Hunter} even proved the sharp lower bound $h_{2p}(a_1, \ldots, a_n) \ge 1/(2^p p!)$ for $a_1^2+\cdots+a_n^2=1$.
Here is an extension of this result to fractional degrees.

\begin{Proposition} \label{abprop1}
Suppose $|\mu-2p| < \tfrac{1}{2}$ for some nonnegative integer $p$ and let $2q$ be the smallest even integer such that
$\mu \le 2q$, i.e., $q=p$ if $\mu \le 2p$ and $q=p+1$ if $\mu > 2p$. Then
\[{\rm Re}\,h_\mu(a_1, \ldots, a_n) \ge \frac{(\mu+n-1)(\mu+n-2)\cdots (\mu+1)}{(2q+n-1)(2q+n-2)\cdots (2q+1)}\frac{\cos(\mu\pi)}{2^q q!}\]
whenever $a_1^2+\cdots+a_n^2=1$.
\end{Proposition}

{\em Proof.}
With $F(x;a_1, \ldots,a_n)$ abbreviated to $F(x)$, we have
\begin{eqnarray*}
\binom{\mu\!+\!n\!-\!1}{n\!-\!1}^{-1}{\rm Re}\,h_\mu(a_1,\ldots,a_n) \!\!
& = & \!\!{\rm Re}\,\left(\int_{-\infty}^0 e^{i\mu\pi}|x|^\mu F(x)\,dx +\int_0^\infty |x|^\mu F(x)\,dx\right)\!\!\\
& = & \!\!\cos(\mu\pi)\int_{-\infty}^0 |x|^\mu F(x)\,dx+\int_0^\infty |x|^\mu F(x)\,dx\!\!\\
& \ge & \! \!\cos(\mu\pi) \int_\mathbb{R}|x|^\mu F(x)\,dx.
\end{eqnarray*}
The equality $a_1^2+\cdots+a_n^2=1$ implies that $|a_j| \le 1$ for all $j$. Thus $[a_1,a_n] \subset [-1,1]$,
and since $|x|^\mu \ge |x|^{2q}$ for $|x| \le 1$, it follows that
\begin{eqnarray*}
\binom{\mu\!+\!n\!-\!1}{n\!-\!1}^{-1}{\rm Re}\,h_\mu(a_1,\ldots,a_n)
& \ge &\cos(\mu\pi)\int_{a_1}^{a_n} |x|^\mu F(x)\,dx\\
& \ge &
\cos(\mu\pi)\int_{a_1}^{a_n} |x|^{2q} F(x)\,dx.
\end{eqnarray*}
But the last integral equals $\tbinom{2q+n-1}{n-1}^{-1}h_{2q}(a_1, \ldots,a_n)$ and Hunter~\cite{Hunter} showed that $h_{2q}(a_1, \ldots,a_n)$ is at least $1/(2^q q!)$.
$\;\:\square$

\medskip
The imaginary part of $h_\mu(a_1,\ldots, a_n)$ is
\[\binom{\mu+n-1}{n-1}\left(\sin(\mu\pi)\int_{-\infty}^0 |x|^\mu F(x)\,dx+\int_0^\infty |x|^\mu F(x)\,dx\right).\]
If $2p < \mu < 2p+1$ with a nonnegative integer $p$, this is strictly positive with the lower bound
\[\frac{(\mu+n-1)(\mu+n-2)\cdots (\mu+1)}{(2p+n+1)(2p+n)\cdots (2p+3)}\frac{\sin(\mu\pi)}{2^q q!}\]
for $a_1^2+\cdots+a_n^2=1$. (Note that the smallest even integer greater than $\mu$ is $2q=2p+2$.)
Thus, if $\mu \in (2p,2p+\tfrac{1}{2})$,
then $h_\mu$ maps all of $\RR^n \setminus \{(0,\ldots,0)\}$ into the open upper-right quarter-plane.
The set $(0,\infty)^n$ is always mapped into the open right half-line. The function $h_\mu$ maps
$(-\infty,0)^n$ into the upper-left quarter-plane for $\mu \in (2p+\tfrac{1}{2},2p+1)$, into the
lower-left quarter-plane for $\mu \in (2p+1,2p+3/2)$, and into the lower-right quarter-plane for $\mu \in (2p+3/2,2p+2)$.

Let again ${\rm Re}\,z > -1$ and let the branch of the complex logarithm be the one specified in Theorem~\ref{abthm2}.
If $\lambda >0$, then $(\lambda a)^z=\lambda^za^z$, but if $\lambda <0$ and $a<0$, then $(\lambda a)^z=\lambda^za^ze^{-2\pi i z}$.
Thus, $h_z(a_1, \ldots,a_n)$ is positively homogeneous but in general not genuinely homogeneous. If $z=\mu$ is a real number and
if $\lambda >0$, we have
\[
{\rm Re}\,h_\mu(\lambda a_1, \ldots, \lambda a_n)={\rm Re}\,[\lambda^\mu h_\mu(a_1, \ldots, a_n)]
=\lambda^\mu\, {\rm Re}\, h_\mu(a_1, \ldots,a_n),
\]
and hence ${\rm Re}\, h_\mu(a_1, \ldots,a_n)$ is also positively homogeneous. This makes Proposition~\ref{abprop1} useful.
However, if, for instance, $z=i\nu$ with a real number $\nu \neq 0$, then, for $\lambda >0$,
\begin{eqnarray*}
& & h_{i\nu}(\lambda a_1, \ldots, \lambda a_n)=\lambda^{i\nu}\,h_{i\nu}(a_1, \ldots, a_n)\\
& & \quad = \Big(\cos(\nu\log \lambda)+i \sin(\nu \log\lambda)\Big)
\Big({\rm Re}\, h_{i\nu}(a_1, \ldots, a_n) +i\, {\rm Im}\,h_{i\nu}(a_1, \ldots, a_n)\Big),
\end{eqnarray*}
which reveals that neither ${\rm Re}\, h_{i\nu}(a_1, \ldots,a_n)$ nor ${\rm Im}\, h_{i\nu}(a_1, \ldots,a_n)$ is positively
homogeneous. The following proposition completes the picture provided by Theorem~\ref{abthm2}.

\begin{Proposition} \label{abprop2}
If $z \in \CC \setminus \RR$ and ${\rm Re}\,z >-1$, then both the real part and the imaginary
part of $h_{z}(a_1, \ldots,a_n)$ are indefinite.
\end{Proposition}

{\em Proof.} From (\ref{aaaa}) we infer that
if $z=\mu+i\nu$ with $\mu,\nu \in \RR$ and $\nu\neq 0$, then, for $a>0$,
\[h_z(a, \ldots, a)=\binom{z+n-1}{n-1}a^{\mu+i\nu}=\binom{z+n-1}{n-1}a^\mu e^{i\nu \log a},\]
which shows that the range of $h_z$ contains a spiral (a circle for $\mu=0$) rotating around the origin and hence reveals that
both ${\rm Re}\,h_z$ and ${\rm Im}\, h_z$ assume strictly positive as well as strictly negative values. $\;\:\square$

\section{Proofs of Theorems \ref{abthm3} and \ref{abthm4}}
\label{bosec5}

\noindent
{\em Proof of Theorem \ref{abthm3}.} If $H(a,a, \ldots,a) \le 0$ for some
$a$ in $(r,s) \setminus\{0\}$, then the inequality $H(a_1, \ldots,a_n)>0$ is not true
for all $(a_1, \ldots, a_n)$ in $(r,s)^n\setminus\{(0,\ldots,0)\}$.

So assume $H(a,a, \ldots, a) >0$ for $a$ in $(r,s) \setminus\{0\}$.
We have to show that then $H(a_1, \ldots,a_n) >0$ whenever $a_j \in (r,s)$ for all $j$ and at least two of
the numbers are different.
Since $H(a_1, \ldots, a_n)$ is symmetric, we may assume
that $a_1 \le \cdots \le a_n$.
We know that Theorem \ref{abthm1} extends to the case~(\ref{aaaaa}).
Thus, we have
\begin{equation} \label{bo51}
H(a_1, \ldots, a_n)=\int_{a_1}^{a_n} \cP(x) F(x;a_1, \ldots, a_n)\,dx
\end{equation}
with
\[\cP(x)=\sum_{j=1}^m \binom{j+n-1}{n-1}c_j x^j.\]
From (\ref{aaaa}) we see that $\cP(x)=H(x,x, \ldots,x)$. Thus, if $H(a,a, \ldots,a) >0$
for $a$ in $(r,s) \setminus\{0\}$, then $\cP(x) >0$ for $x \in (r,s)\setminus\{0\}$
and~(\ref{bo51}) implies that $H(a_1, \ldots,a_n) >0$ if $r < a_1 \le \cdots \le a_n < s$
and at least two of the $a_j$ are different. $\;\:\square$

\medskip
\noindent
{\em Proof of Theorem \ref{abthm4}.} Since we require that $H(a,a, \ldots, a) >0$ for nonzero $a \in (r,s)$,
we are left with the case where $r < a_1 \le \cdots \le a_n < s$ and $a_1 < a_n$. We then get that $H(a_1, \ldots,a_n)$
equals
\[\int_{a_1}^{a_n}\int_{a_1}^{a_n} \cP(x,y) F(x;a_1, \ldots,a_n) F(y; a_1, \ldots, a_n)\, dx \,dy\]
with
\[\cP(x,y)=\sum_{j,k=1}^m \binom{j+n-1}{n-1}\binom{k+n-1}{n-1}c_{jk} x^jy^k.\]
From (\ref{aaaa}) it follows that
\[\cP(x,y)=\sum_{j,k=1}^m c_{jk} h_j(x,x, \ldots,x) h_k(y,y, \ldots,y).\]
Consequently, if $\cP(x,y) \ge 0$ on $(r,s)^2 \setminus \{(0,0)\}$, then the double integral is strictly positive. $\;\:\square$

\section{Emergence of Schur Polynomials}
\label{absec4}

\noindent
Throughout the following think of $a_1, \ldots,a_n$ as variables or as nonzero and pairwise distinct real numbers.
Given an $n$-tuple $\la=(\la_1, \la_2, \ldots,\la_n)$ of integers satisfying $\la_1 \ge \la_2 \ge \cdots \ge \la_n \ge 0$,
the {\em Schur polynomial} $s_\la(a_1, a_2, \ldots,a_n)$ is defined as
\begin{equation}\label{abSchur}
s_\la(a_1, a_2,\ldots,a_n)\\
= \dfrac{\det \begin{bmatrix} a_1^{\lambda_n} & a_1^{\lambda_{n-1}+1} & a_1^{\lambda_{n-2}+2} & \cdots & a_1^{\lambda_1+n-1} \\ a_2^{\lambda_n} & a_2^{\lambda_{n-1}+1} & a_2^{\lambda_{n-2}+2} & \cdots & a_2^{\lambda_1+n-1}  \\ \vdots & \vdots & \vdots & \ddots & \vdots \\ a_n^{\lambda_n} & a_n^{\lambda_{n-1}+1} & a_n^{\lambda_{n-2}+2} & \cdots & a_n^{\lambda_1+n-1} \end{bmatrix}}{ V(a_1,a_2,\ldots,a_n)};
\end{equation}
see, for example, \cite{ec2}.
From (\ref{abJac}) we see that if $z$ is a nonnegative integer, then
\[h_z(a_1,a_2,\ldots,a_n)=s_{(z,0,\ldots,0)}(a_1,a_1, \ldots,a_n),\]
with $s_{(0,0,\ldots,0)}(a_1, a_2,\ldots,a_n)=1$.

\begin{Proposition} \label{Ex1}
Let $z$ be a positive integer. If $1 \le z \le n-1$, then $h_{-z}(a_1, \ldots,a_n)=0$. If $z \ge n$,
then
\[h_{-z}(a_1, \ldots,a_n)=(-1)^{n-1}(a_1\cdots a_n)^{n-1-z}s_{(z-n, \ldots, z-n,0)}(a_1, \ldots,a_n).\]
\end{Proposition}

{\em Proof.} Consider (\ref{abJac}) with $z$ replaced by $-z$. If $1 \le z \le n-1$, then the determinant on the right
contains a repeated column and hence it is zero. So let $z \ge n$. Then, again by~(\ref{abJac}),
\[h_{-z}(a_1,\ldots,a_n)  V(a_1,\ldots,a_n) = \det
\begin{bmatrix}
1 & a_1 & a_1^2 & \cdots & a_1^{n-2} & a_1^{-z+n-1} \\
1 & a_2 & a_2^2 & \cdots & a_2^{n-2} & a_2^{-z+n-1} \\
\vdots & \vdots & \vdots & \ddots & \vdots & \vdots \\
1 & a_n & a_n^2 & \cdots & a_n^{n-2} & a_n^{-z+n-1} \\
\end{bmatrix},\]
and this equals $(a_1\cdots a_n)^{-z+n-1}$ times
\[ \det
\begin{bmatrix}
a_1^{0+(1+z-n)} & a_1^{1+(1+z-n)} & a_1^{2+(1+z-n)} & \cdots & a_1^{n-2+(1+z-n)} & 1 \\
a_2^{0+(1+z-n)} & a_2^{1+(1+z-n)} & a_2^{2+(1+z-n)} &  \cdots & a_2^{n-2+(1+z-n)} & 1 \\
\vdots & \vdots  & \ddots & \vdots & \vdots \\
a_n^{0+(1+z-n)} & a_n^{1+(1+z-n)} & a_n^{2+(1+z-n)} & \cdots & a_n^{n-2+(1+z-n)} & 1 \\
\end{bmatrix}.
\]
This last determinant is
\[
(-1)^{n-1}\det
\begin{bmatrix}
a_1^0 & a_1^{1+(z-n)} & a_1^{2+(z-n)} & \cdots & a_1^{(n-1)+(z-n)} \\
a_2^0 & a_2^{1+(z-n)} & a_2^{2+(z-n)} & \cdots & a_2^{(n-1)+(z-n)} \\
\vdots & \vdots & \vdots & \ddots & \vdots \\
a_n^0 & a_n^{1+(z-n)} & a_n^{2+(z-n)} & \cdots & a_n^{(n-1)+(z-n)} \\
\end{bmatrix}.
\]
Thus, letting \[\lambda=(\underbrace{z-n,z-n,\ldots,z-n}_{n-1 \mbox{ copies}},0)\] we get
\[
h_{-z}(a_1,\ldots,a_n)=(-1)^{n-1}(a_1 \cdots a_n)^{n-1-z} s_{\lambda}(a_1,\ldots,a_n).\;\:\square
\]

\begin{Proposition} \label{Ex2}
Let $z$ be a positive rational number but not be an integer. Write $z=p/q$ with $q \ge 2$ and $gcd(p,q)=1$. Then
$h_z(a_1,\ldots,a_n)$ is
\[\prod_{1 \leq i < j \leq n} \frac{1}{a_i^{(q-1)/q}+a_i^{(q-2)/q}a_j^{1/q}+ \cdots + a_j^{(q-1)/q}}  s_{\lambda}(a_1^{1/q},\ldots,a_n ^{1/q}).\]
\end{Proposition}

{\em Proof.}
We start again with (\ref{abJac}). The determinant on the right may be written as
\[\det
\begin{bmatrix}
1 & (a_1^{1/q})^q & (a_1^{1/q})^{2q} & \cdots & (a_1^{1/q})^{(n-2)q} & (a_1^{1/q})^{p+(n-1)q} \\
1 & (a_2^{1/q})^q & (a_2^{1/q})^{2q} & \cdots & (a_2^{1/q})^{(n-2)q} & (a_2^{1/q})^{p+(n-1)q} \\
\vdots & \vdots & \vdots & \ddots & \vdots & \vdots \\
1 & (a_n^{1/q})^q & (a_n^{1/q})^{2q} & \cdots & (a_n^{1/q})^{(n-2)q} & (a_n^{1/q})^{p+(n-1)q} \\
\end{bmatrix}.\]
This equals
\begin{eqnarray*}
& & \det \left[\begin{array}{cccc}
1 & (a_1^{1/q})^{1+(q-1)} & (a_1^{1/q})^{2+2(q-1)} & \cdots  \\
1 & (a_2^{1/q})^{1+(q-1)} & (a_2^{1/q})^{2+2(q-1)} & \cdots  \\
\vdots & \vdots & \vdots & \ddots  \\
1 & (a_n^{1/q})^{1+(q-1)} & (a_n^{1/q})^{2+2(q-1)} & \cdots  \\
\end{array}\right.\\
& & \left. \begin{array}{ccc}
\qquad \qquad\cdots & (a_1^{1/q})^{n-2+(n-2)(q-1)} & (a_1^{1/q})^{(n-1)+p+(n-1)(q-1)} \\
\qquad \qquad \cdots & (a_2^{1/q})^{n-2+(n-2)(q-1)} & (a_2^{1/q})^{(n-1)+p+(n-1)(q-1)} \\
\qquad \qquad \ddots & \vdots & \vdots \\
\qquad \qquad \cdots & (a_n^{1/q})^{n-2+(n-2)(q-1)} & (a_n^{1/q})^{(n-1)+p+(n-1)(q-1)} \\
\end{array}\right],
\end{eqnarray*}
and from (\ref{abSchur}) we deduce that the last determinant is
\[V(a_1^{1/q},\ldots,a_n^{1/q})  s_{\lambda}(a_1^{1/q},\ldots,a_n^{1/q})\]
with $\la=(p+(n-1)(q-1),(n-2)(q-1),\ldots,2(q-1),(q-1),0)$. Consequently,
\begin{eqnarray*}
& & h_z(a_1,\ldots,a_n) =\frac{\det V(a_1^{1/q},\ldots,a_n^{1/q})}{\det V(a_1,\ldots,a_n)} s_{\lambda}(a_1^{1/q},\ldots,a_n^{1/q}) \\
& &= \prod_{1 \leq i < j \leq n} \frac{1}{a_i^{(q-1)/q}+a_i^{(q-2)/q}a_j^{1/q}+ \cdots + a_j^{(q-1)/q}} \cdot s_{\lambda}(a_1^{1/q},\ldots,a_n ^{1/q}). \;\:\square
\end{eqnarray*}

These ideas extend to a related formula when $p/q$ is negative.
We leave the details to the interested reader.

\begin{Example} \label{Ex3}
{\rm If $z=2/3$ and $n=4$,  then $\lambda=(2+3\cdot 2, 2 \cdot 2, 2, 0)=(8,4,2,0)$ and we obtain that}
\begin{align*}
&h_{\frac{2}{3}}(a_1,a_2,a_3,a_4)
\\ &= \left( \prod_{1 \leq i < j \leq 4} \frac{1}{a_i^{2/3}+a_i^{1/3}a_j^{1/3}+a_j^{2/3}} \right) \cdot s_{(8,4,2,0)}(a_1^{1/3},a_2^{1/3},a_3^{1/3},a_4^{1/3}).
\end{align*}
\end{Example}

\medskip
{\bf Acknowledgment.} We thank Grigori Olshanski and Terence Tao 
for their valuable comments. In particular, Grigori Olshanski's hint to~\cite{CS} solved the problem (in the affirmative)
whether $F(x;a_1, \ldots,a_n)$ is unimodal, which was left as an open question in~\cite{GOO}.

\end{document}